\documentclass[11 pt, reqno]{article}

\usepackage{amsmath,amsfonts,amssymb,amsthm}
\usepackage[colorlinks=true,hyperindex=true]{hyperref}

\usepackage[a4paper, left=2cm, right=2cm, top=2 cm, bottom= 2 cm]{geometry}
\usepackage{color}
\usepackage{fancyhdr}
\usepackage{latexsym}

\usepackage{chngcntr}
\counterwithin*{equation}{section}

\newtheorem{ccounter}{ccounter}[section]
\newtheorem{thm}[ccounter]{Theorem}
\newtheorem{lem}[ccounter]{Lemma}
\newtheorem{cor}[ccounter]{Corollary}
\newtheorem{defn}[ccounter]{Definition}
\newtheorem{prop}[ccounter]{Proposition}
\newtheorem{ass}[ccounter]{Assumption}
\newtheorem{ex}[ccounter]{Example}

\def\bet{\begin{thm}}
\def\eet{\end{thm}}
\def\bel{\begin{lem}}
\def\eel{\end{lem}}
\def\bas{\begin{ass}}
\def\eas{\end{ass}}
\def\bec{\begin{cor}}
\def\eec{\end{cor}}
\def\bed{\begin{defn}}
\def\eed{\end{defn}}
\def\bep{\begin{prop}}
\def\eep{\end{prop}}
\def\beq{\begin{equation}}
\def\eeq{\end{equation}}
\def\proof{\noindent {\bf Proof.}\ \ }
\def\bea{\begin{equation*}}
\def\eea{\end{equation*}}
\def\tr{\mathrm{tr}}
\def\bex{\begin{ex}}
\def\eex{\end{ex}}
\def\remark{\noindent{\bf Remark. }}
\def\bfx{\textbf{x}}
\def\bfy{\textbf{y}}
\def\bfz{\textbf{z}}
\def\mfa{\mathfrak{a}}
\def\Oma{\Omega_{\mathfrak{a}}}

\def\rr{\mathbb{R}}

\def\cc{\mathbb{C}}
\def\1{\boldsymbol{1}}
\def\Im{\mathrm{Im}}
\def\Re{\mathrm{Re}}
\def\e{\mathrm{e}}
\def\i{\mathrm{i}}
\def\del{\partial}
\def\d{\mathrm{d}}
\def\eps{\varepsilon}
\renewcommand\leq\varleq
\renewcommand\geq\vargeq
\def\ee{\mathrm{E}}

\def\O{\mathcal{O}}

\def\ee{\mathbb{E}}
\def\om{\omega}
\def\pp{\mathbb{P}}

\def\tilf{\tilde{f}}
\def\convdist{\overset{d}{\to}}
\def\rhosc{\rho_{\mathrm{sc}}}
\def\gV{\gamma^{(V)}}

\def\gU{\gamma^{(U)}}
\def\rhoU{\rho^{(U)}}
\def\rhoV{\rho^{(V)}}
\def\hatp{\hat{p}}
\def\gf{\gamma^{(\mathfrak{f})}}
\def\hatzeta{\hat{\zeta}}
\def\bx{\boldsymbol{x}}
\def\by{\boldsymbol{y}}
\def\ea{e_{\mathfrak{a}}}
\def\Ea{E_{\mathfrak{a}}}
\def\I{\mathcal{I}}
\def\num{\mathcal{N}}
\def\Lone{1}

\def\G{\Phi}

\def\A{\mathcal{A}}

\begin{document}
\title{ Applications of mesoscopic CLTs in Random Matrix Theory }


\begin{table}
\centering

\begin{tabular}{c}
\multicolumn{1}{c}{\Large{\bf Applications of mesoscopic CLTs in Random Matrix Theory }}\\
\\
\\
\end{tabular}
\begin{tabular}{c c c}
Benjamin Landon & &  Philippe Sosoe\\
 & & \\
{\small Massachusetts Institute of Technology }& & {\small Cornell University}\\
& & \\
\end{tabular}
\\
\begin{tabular}{c}
\multicolumn{1}{c}{\today}\\
\\
\end{tabular}

\begin{tabular}{p{15 cm}}
\small{{\bf Abstract:} We present some applications of central limit theorems on mesoscopic scales for random matrices.  When combined with the recent theory of ``homogenization'' for Dyson Brownian Motion, this yields the universality of quantities which depend on the behavior of single eigenvalues of Wigner matrices and $\beta$-ensembles.  Among the results we obtain are the Gaussian  fluctuations of single eigenvalues for Wigner matrices (without an assumption of 4 matching moments) and classical $\beta$-ensembles ($\beta=1, 2, 4$), 
Gaussian fluctuations of the eigenvalue counting function, and an asymptotic expansion up to order $o(N^{-1})$ for the expected value of eigenvalues in the bulk of the spectrum.  The latter result solves a conjecture of Tao and Vu.}
\end{tabular}
\end{table}


\section{Introduction} The {\let\thefootnote\relax\footnote{P.S. is partially supported by NSF grant DMS-1811093.}object of this note is to show how two recent results in random matrix theory, central limit theorems for linear statistics of eigenvalues on mesoscopic scales (see \cite{heknowles, bekermanlodhia}), and homogenization for Dyson Brownian Motion (DBM) introduced in \cite{fixed1} and refined in \cite{fixed2}, can be used to derive a number of results concerning the fluctuations of single eigenvalues. In particular, we obtain the following results:
\begin{enumerate}
\item Gaussian fluctuations on scale $\frac{\sqrt{\log N}}{N}$ for individual eigenvalues of Wigner matrices in the bulk,
\item Gaussian fluctuation on scale $\frac{\sqrt{\log N}}{N}$ for individual eigenvalues of classical $\beta$-ensembles in the bulk ($\beta=1, 2, 4$),
\item Gaussian fluctuations for eigenvalue counting functions,
\item Asympotic expansion of the mean and variance of eigenvalues in the bulk, up to an error of order $o(N^{-1})$.
\end{enumerate}
Here, $N$ denotes the dimension of the random matrix ensemble under consideration. 
Item 1 was first obtained for the Gaussian Unitary Ensemble (GUE) by Gustavsson \cite{gustavsson}, and then for complex Hermitian Wigner matrices whose moments match the GUE to fourth order by Tao and Vu \cite{TV}.  O'Rourke \cite{orourke} extended this to the Gaussian Orthogonal and Symplectic Ensembles (GOE/GSE) by applying a result of Forrester and Rains \cite{FR} (and then to Wigner matrices of these symmetry classes with four matching moments).   
Gustaffson originally first proved Item 3 and then obtained Item 1  via the duality
$$
| \{ i : \lambda_i < y \} | \ge k \iff \lambda_k < y,
$$
where the $\{ \lambda_i \}_i$ are random matrix eigenvalues in increasing order. 
Dallaporta and Vu \cite{dallaportavu} used this to obtain the central limit theorem for the eigenvalue counting function of Wigner matrices (of those with four moments matching the Gaussian ensembles) and additionally computed the asymptotic expectation and variance of this quantity.

In the present work we remove all moment matching assumptions in Items 1 and 3.  Moreover our approach also applies to $\beta$-ensembles and we obtain the Gaussian fluctuations for the eigenvalues of the classical $\beta$-ensembles.  Our proof is based on \emph{universality}, meaning that we show that the eigenvalue fluctuations are universal, coinciding with the Gaussian case.  In the $\beta$-ensemble case, this approach applies for even non-classical values of $\beta \ge 1$, and so if the analog of Gustavsson's result were known for general Gaussian $\beta$-ensembles, then we would obtain Item $2$ for all $\beta \ge 1$.
 
 The formula we obtain for the expected value of a single bulk eigenvalue in Item 4 solves a conjecture of Tao and Vu \cite[Conjecture 1.7]{taovu}, and extends it to the real symmetric case.

\subsection{Homogenization of Dyson Brownian Motion}

Our work relies on a recent technical innovation in the study of Dyson Brownian motion (DBM), known as ``homogenization.''  DBM is a stochastic process on random matrices first introduced by Dyson \cite{D}.  It was first applied by Erd\H{o}s, Schlein and Yau \cite{ESY} who introduced DBM as a tool to study local eigenvalue fluctuations and prove Wigner-Dyson-Mehta conjecture for all symmetry classes.  

  Homogenization relies on a coupling between 
 two Dyson Brownian motions $\bx(t)$ and $\by (t)$.  The first process $\bx (t)$ has initial data given by the Wigner matrix whose eigenvalues we wish to study, and $\by (t)$ comes from the equilibrium Gaussian ensemble that we would like to compare $\bx (t)$ to. 
 The difference between these two processes satisfies a discrete parabolic equation; the work \cite{fixed1} of Bourgade, Erd\H{o}s, Yau and Yin establishes a homogenization theory for this parabolic equation, resulting in the estimate
\begin{equation}\label{eqn: rep}
x_i(t)-y_i(t)=\frac{1}{N}\zeta_{x}-\frac{1}{N}\zeta_{y}+\O(N^{-1-c}).
\end{equation}
for some $c>0$ with probability $1 - o(1)$. 
Here,
\begin{equation}\label{eqn: Gls}
\zeta_{x}=\sum_{j=1}^N \G(x_j(0))-\sum_{j=1}^N \G (\gamma_j), \qquad \zeta_{y}=\sum_{j=1}^N \G(y_j(0))-\sum_{j=1}^N \G(\gamma_j)
\end{equation}
for some explicit, $N$-dependent function $\G$.   Above, the $\gamma_j$ are the quantiles of the semicircle density, the limit of the empirical eigenvalue density of $\bx (0)$ and $\by (0)$, and so the second sums are roughly the expected value of the first sums.  The function $G$ is constructed from the heat kernel of a nonlocal operator that arises as a continuum approximation to the equation satisfied by $\bx-\by$ (hence the name ``homogenization''). 

The linear statistic \eqref{eqn: Gls} is on the  scale $t$ which must satisfy $1 \gg t \gg N^{-1}$ for the estimate \eqref{eqn: rep} to hold.  This scale is in between the microscopic and global scales, and so the quantity \eqref{eqn: Gls} is known as a mesoscopic linear statistic.

The representation \eqref{eqn: rep} allows for a precise description of the fluctuations of quantities depending on single eigenvalues $x_i(t)$, provided one has good control of mesoscopic linear statistics, specifically, the quantities  \eqref{eqn: Gls}. 
At the time of the publication of \cite{fixed1}, such a general result was not yet available, and the authors found that it sufficed, for their intended application, to obtain a result about mesoscopic linear statistics for the GOE.  The intended application of homogenization in \cite{fixed1} was to the local correlation functions, which concern the scale $N^{-1}$ and differs crucially from the $\sqrt{ \log(N)} N^{-1}$ of the single eigenvalue fluctuations.  For this reason, it was not necessary to treat the quantities \eqref{eqn: Gls} for the non-GOE case.

  In \cite{fixed2}, the authors of the present work (with Yau) revisited the argument in \cite{fixed1} and derived a refined homogenization result depending only on local assumptions on the initial data.    We also obtained precise control on linear statistics of observables such as $\G$ in \eqref{eqn: Gls} in the case of deformed Wigner matrices.   The derivation of our main results will follow from a combination of \eqref{eqn: rep} and a strengthened version 
  of the mesocopic CLT we proved in \cite{fixed2}.

\subsection{Mesoscopic linear statistics}
Let $\lambda_1 \le \lambda_2 \cdots \le \lambda_N$, denote the eigenvalues of a 
 Wigner matrix. It is by now classical that for smooth functions $f$, the quantity
\[\sum_{j=1}^N f(\lambda_j)-\mathbb{E}\left[\sum_{j=1}^N f(\lambda_j) \right] \]
is asymptotically Gaussian for smooth functions $f$. A general result applying to sufficiently regular functions 
 appears in the work of Lytova and Pastur \cite{lytovapastur}.
For compactly supported or rapidly decaying functions $f$, $0<\alpha<1$ and an energy $E\in (-2,2)$,
\begin{equation}\label{eqn: sums}
  \sum_{j=1}^N f(N^\alpha (\lambda_j-E))-\mathbb{E}\left[\sum_{j=1}^N f(N^\alpha (\lambda_j-E)) \right]
\end{equation}
is asymptotically Gaussian, with limiting variance
\beq \label{eqn:int_var}
\frac{c_{sym}}{2\pi^2}\int\int \left(\frac{f(x)-f(y)}{x-y}\right)^2\,\mathrm{d}x\mathrm{d}y.
\eeq
Here $c_{sym}=\frac{1}{2}$ for Hermitian Wigner matrices, and $c_{sym}=1$ for real symmetric Wigner matrices. Because of the rescaling by $N^\alpha$, the sums \eqref{eqn: sums} typically involve about $N^{1-\alpha}$ eigenvalues, and consequently such linear statistics are refered to as \emph{mesoscopic}. This result has been known in restricted cases for some time (see \cite{ak1,ak2,lodhiasimm}).  In \cite{heknowles}, He and Knowles obtain the result for general Wigner matrices and general $f$. See also the work of Bekerman and Lodhia \cite{bekermanlodhia} for the case of $\beta$-ensembles. We gave  proofs of  more general results in \cite{fixed2} in the deformed Wigner case, and one-cut $\beta$-ensembles.

The method we use to derive our mesoscopic CLT represents a considerable simplification over the proof offered in \cite{fixed2}.  Moreover, the result we derive here is stronger than that in \cite{fixed2}, which had some restrictions on the scale of the function under consideration - this was an artifact of the proof given there, and absent from the present work.   

This strengthening is in  large part due to the fact that we are dealing with Wigner matrices, rather than the deformed Wigner matrices of our previous paper.   
Our result improves on the known results for Wigner matrices \cite{heknowles,bekermanlodhia}  in that we allow the function $f$ to not have compact support or decay quickly at infinity. As a consequence, the variance of the linear statistic may grow logarithmically with $N$. This is essential for our application to \eqref{eqn: Gls} - the function $G$ turns out to roughly behave as a smoothed out step function on the mesoscopic scale $t$, and for such a function the quantity \eqref{eqn:int_var} is $\O ( \log (N) )$.  Note that the size of the fluctuations of the linear statistics \eqref{eqn: Gls} and the single eigenvalues are the same, and so the mesoscopic central limit theorem is necessary for treating the single eigenvalue fluctuations via the homogenization approach.   Such a CLT was not known prior to the present work. 

The proof of the CLT, Proposition \ref{prop:stein}, is inspired by the method of Shcherbina \cite{shcherbina}, but uses the cumulant expansion (as in \cite{heknowles}) instead of the more intricate resolvent expansion in \cite{fixed2}.  We also rely on the isotropic local law of \cite{BEKYY,knowlesyin} to simplify the treatment of some error terms (this is not strictly necessary as such a result was not available for \cite{fixed2}).

\subsection{Statement of results}
For simplicity, we state all our results for Wigner matrices in the case when $H$ is real and symmetric. The same proofs apply with only minor modification to random Hermitian matrices.  Where relevant, we comment on the corresponding results for the complex Hermitian case.   In this case however, many of results below can be derived using the Br\'ezin-Hikami formula \cite{BH}, requiring neither of the techniques highlighted in the introduction.

Let $\xi_o$ and $\xi_d$ be two real centered random variables with bounded moments of all orders, with the variance of $\xi_o$ being $1$.   A real symmetric Wigner matrix is an $N \times N$ self adjoint matrix so that the entries $\{ H_{ij} \}_{i \le j}$ are independent, and 
$$
\sqrt{N} H_{ii} \sim \xi_d, \qquad \sqrt{N} H_{ij} \sim \xi_o, \quad i \neq j.
$$
We also use the following notation for the cumulants of the matrix entries,
\begin{align}
s_k:=\frac{1}{i^k}\frac{\mathrm{d}^k}{\mathrm{d}^kt}\log \mathbb{E}[e^{it \xi_o}] \bigg\vert_{t=0}  \notag \\
s_k+a_k:=\frac{1}{i^k}\frac{\mathrm{d}^k}{\mathrm{d}^kt}\log \mathbb{E}[e^{it \xi_d}] \bigg\vert_{t=0} \label{eqn:cumulants}
\end{align}
Our notational choice for the second line of \eqref{eqn:cumulants} is so that the cumulant expansion we use to derive our mesoscopic central limit theorem has a simple form (see the proof of Lemma \ref{lem:cumulant}).

We denote by $\lambda_1\le \cdots \le \lambda_N$ be the ordered eigenvalues of $H$. As a consequence of the local semi-circle law, Theorem \ref{thm: sc-law}, $\lambda_i$ is typically very close to the $i$th $N$-quantile $\gamma_i$ of the semicircle distribution, also known as the ``classical location'' of $\lambda_i$.  The semicircle and classical eigenvalue locations are defined by
\begin{equation}\label{eqn: classical-loc}
\rho_{\mathrm{sc}}(x)=\frac{1}{2\pi}\sqrt{(4-x^2)_+} , \qquad \int^{\gamma_i}_{-2}\rho_{\mathrm{sc}}(x)\,\mathrm{d}x=\frac{i}{N}.
\end{equation}
Our first result describes the fluctuations of $\lambda_i$ about $\gamma_i$.
\bet \label{thm1}
Let $\kappa>0$, and let $\lambda_i$ be the eigenvalues of a symmetric Wigner matrix.  For any sequence $i = i_N$ such that $\kappa N \le i \le (1- \kappa)N$ we have,
\beq
N\frac{ \lambda_i - \gamma_i}{\sqrt{ \frac{\log(N)}{1-\gamma_i^2/4}}} \convdist N(0, 1).
\eeq
In the case of complex Hermitian Wigner matrices, the same result holds with limiting variance $\frac{1}{2}$.
\eet
As explained in the introduction, this has been obtained only for the Gaussian ensembles and those Wigner matrices whose first four moments match the Gaussian ensembles \cite{orourke,gustavsson,TV}. 
Our contribution is thus to remove the assumption of matching moments. Theorem \ref{thm1} has also been  obtained independently by Bourgade and Krishnan \cite{BK} by different methods.

The result of Gustavsson holds for any sequence of indices $i_N$ such that $i_N \to \infty $ and $|N-i_N| \to \infty$, i.e., for indices close to the edge.  We expect that our universality result holds for any sequences such that $\kappa_N N \le i_N \le (1- \kappa_N) N$ with $\kappa_N \ge N^{-c}$ for some sufficiently small $c>0$.  We do not pursue this here in the interest of brevity. Results for the case of moderately growing $i_N$ follow from edge universality.  For example, Gaussian fluctuations in the case  that $i_N \to \infty$ but $i_N \le N^{c}$ for some sufficiently small $c>0$, follow from \cite{edgebeta}.  

By a well-known argument (see \cite{gustavsson,dallaportavu} or below), the result on fluctuations of eigenvalues implies a CLT for indicator functions:
\bec \label{cor:counting}
For $E \in (-2, 2)$ let $\num(E) = \# \{ i  : \lambda_i \le E \}$ be the eigenvalue counting function of a real symmetric Wigner matrix.  Then,
\beq
\frac{\num(E) - N \int_{-\infty}^E \d \rho_{\mathrm{sc}} }{\pi^{-1} \sqrt{ \log(N) } } \convdist N(0, 1).
\eeq
\eec

A similar method gives the following for classical $\beta$-ensembles in the one-cut case.  We recall that a $\beta$-ensemble with potential $V : \rr \to \rr$ is the measure $\mu_V$ on the simplex $\{ \lambda_1 < \lambda_2 < \dots < \lambda_N \subseteq \rr^N$ proportional to the Gibbs weight with Hamiltonian $\mathcal{H}$, 
$$
\mu_V \to \e^{ - \beta N \mathcal{H}}, \qquad \mathcal{H} := \frac{1}{2} \sum_{j =1}^N V ( \lambda_j ) - \frac{1}{N} \sum_{i < j } \log ( \lambda_j - \lambda_i )
$$
where $V$ is assumed to satisfy a growth condition so that this measure is finite (see, e.q., \cite{AGZ}).  For sufficiently regular $V$, the empirical measure is known to converge a.s. to a determininistic probability measure we call the \emph{equilibrium measure} and denote by $\rhoV$ \cite{AGZ, KM}.  The one-cut case is the case that this measure is supported on a single interval.  
\bet \label{thm: beta}
Let $V$ be a one-cut $\beta$-ensemble with equilibrium measure supported in $[A, B]$ and $V \in C^4$.   Assume that $V$ is regular in the sense of \cite{KM}.  Denote the particles by $\lambda_i$.  Let $i=i_N$ be a sequence of indices satisfying $\kappa N \le i \le (1- \kappa )N$ for fixed $\kappa >0$.  Then, for $\beta=1, 2, 4$ we have

\beq
 \frac{ \beta^{1/2} N \rhoV ( \gV_i ) \pi}{ 2 \sqrt{ \log(N)}} \left( \lambda_ i - \gV_i \right) \convdist N(0,1 ).
\eeq
where $\rhoV$ is the equilibrium measure of $V$ and $\gV_i$ is the $i$th classical eigenvalue location.
\eet
If the result for $\beta \notin \{1, 2, 4\}$ were known for the Gaussian $\beta$-ensemble, then the above result would be true for all $\beta \ge 1$.  Instead, we have the statement that for all $\beta \ge 1$ and smooth test functions $F$,
\beq
\lim_{ N \to \infty} \left| \ee_V [ F (N \rhoV ( \gV_i) ( \lambda_i - \gV_i)/ \sqrt{\log(N)})] - \ee_{G\beta E} [ F (N \rhosc ( \gamma_i ) ( \lambda_i - \gamma_i)/ \sqrt{\log(N)})]  \right| = 0.
\eeq
The expectation on the LHS is with respect to the $\beta$-ensemble with potential $V$ and on the RHS we have the expectation with respect to the Gaussian $\beta$-ensemble.

As in Corollary \ref{cor:counting} we also get a CLT for the eigenvalue counting function for $\beta$ ensembles with classical values of $\beta$.  A similar result for the indicator function $\1_{(a, b)}$ where $(a, b) \subseteq (A+\kappa, B-\kappa)$ was obtained by Shcherbina \cite{Sh1} for classical $\beta$, using Fredholm determinants.

Our method also allows for an asymptotic expansion of $\ee[ \lambda_i]$ for Wigner matrices up to order $o(N^{-1})$, resolving a conjecture of Tao and Vu \cite{taovu}.  A related result for Gaussian divisible ensembles in the complex Hermitian case was obtained by Edelman-Guionnet-P\'ech\'e \cite{EGP} using the Br\'ezin-Hikami formula. This result finds an expansion for the quantiles of the expected density of states; however this result does not imply the corresponding expansion for the expected eigenvalue location.
\bet \label{thm:evalues} Fix $\kappa >0$.  Let $H$ be a real symmetric Wigner matrix as above, and $\kappa N \le i \le (1-\kappa)N$.  Then,
\beq \label{eqn:expect}
N \ee_H [ \lambda_i-\gamma_i ] =\frac{1}{2 \pi \rhosc ( \gamma_i ) } \arcsin \left( \frac{ \gamma_i}{2} \right) - \frac{1}{ 2 \rhosc ( \gamma_{i} )}
+ \frac{s_4}{4 } ( \gamma_i^3-2 \gamma_i) + \frac{a_2-1}{2} \gamma_i + o (1)
\eeq
and
\beq \label{eqn:variance}
\mathrm{Var}_H ( \lambda_i ) = \mathrm{Var}_{GOE} (\lambda_i ) +\frac{s_4}{8 N^2}\gamma_i^2 +\frac{a_2-1}{N^2}+ o (N^{-2} ).
\eeq
\eet
\remark Tao and Vu \cite{taovu} conjectured in the complex Hermitian case a formula similar to \eqref{eqn:expect}. The formula we find sharpens their prediction, by identifying the quantity $C_{i,N}$ in their proposed asympotic $\ee_{GUE} [ \lambda_i ] = \gamma_i + N^{-1} C_{i, N} + o (N^{-1})$. Their main interest was the dependence of the expectation on the fourth moment of the entries, and we find exactly the dependence predicted in \cite{taovu}.   In the complex Hermitian analog of \eqref{eqn:expect}, the $\arcsin ( \gamma_i )$ term is not present, and $s_4$ is replaced by the sum of the fourth cumulants of the real and imaginary parts of the off-diagonal entries of the matrix.  The formula \eqref{eqn:variance} appears new.  \qed

For a $1$-cut $\beta$-ensemble with equilibrium distribution $\rho_V$ supported on $[A, B]$ we obtain the following:
\bet \label{thm:betaexpect1} Fix $ \kappa >0$. Let $i$ be an index so that $ \kappa N \le i \le (1- \kappa )N$.  We have,
\begin{align}
N\ee_V [ \rho_V ( \gV_i ) ( \lambda_i  - \gV_i) ]
=& - \frac{1}{2 }+ \frac{1}{2\pi^2}\left(\frac{2}{\beta}-1\right)\int_A^{\gV_i} \d \nu_V (x) 
\end{align}
where $\nu_V (x)$ is a specific signed measure, described in more detail below.
\eet

The mesoscopic central limit theorem we present, Theorem \ref{thm:main-meso-clt}, is fairly robust. To illustrate this, we state an extension of some previous results on ``partial linear statistics'' of Bao, Pan and Zhou \cite{bpz}. Specifically, we weaken the regularity assumptions on the functions $f$ appearing there, and do not require the moment matching hypotheses of \cite{bpz}.
\bet \label{thm: BPZ}
Let $H$ be a real symmetric Wigner matrix.  Fix $u \in (-2, 2)$.  Let $f$ be a $C^3$ function such that $f'(x) \neq 0$ only for $x \in (-2+ \kappa, 2 - \kappa )$, and $f(u) = 0$.  Then,
\beq\label{eqn: BPZ-1}
\sum_i f ( \lambda_i ) \1_{ \{ \lambda_i \le u \} } - \ee[ \sum_i f ( \lambda_i ) \1_{ \{ \lambda_i \le u \} }] \convdist N(0, \sigma^2)
\eeq
where $\sigma^2$ is an explicit function of $f$. In fact,
\[\sigma^2= V(f( \cdot )\mathbf{1}_{\ \cdot\  \le u}),\]
where $V$ is defined in \eqref{eqn:EV}.

Suppose $k/N\rightarrow u \in (-2,2)$. Then
\beq\label{eqn: BPZ-2}
\sum_{i=1}^k f ( \lambda_i) - \ee[ \sum_{i=1}^k f ( \lambda_i)] \to N(0, \sigma'^2),
\eeq
where $\sigma'^2= V(f(\cdot)\mathbf{1}_{ \ \cdot\ \le \tilde{\gamma}_u})$, where the function $\tilde{\gamma}_u$ is defined by
\[\int_{-2}^{\tilde{\gamma}_u}\rho_{\mathrm{sc}}(x)\,\mathrm{d}x=u.\]
\eet
\remark If $f(u)\neq 0$, the quantity \eqref{eqn: BPZ-1} has variance of order $\log N$; a central limit theorem for this case follows from our result for the eigenvalue counting function.  Further, we introduced the notation $\tilde{\gamma}_u$ as the inverse of the cumulative distribution function of the semicircle distribution so as to not conflict with our earlier notation for the $N$-quantiles $\gamma_k$.  Note,
$$
\gamma_k = \tilde{\gamma}_{k/N}.
$$

\noindent{\bf Acknowledgements.} We thank the anonymous referee for a careful reading of this work and numerous comments improving the presentation.

\section{The local semi-circle law} \label{sec:lsc}
Let $H$ be a Wigner matrix as above. The empirical distribution of the eigenvalues of such a matrix converges to the semi-circle distribution \eqref{eqn: classical-loc}. 
The \emph{local semi-circle law}, which will use throughout without further comment represents a considerable strengthening of this statement. It is expressed in terms of the \emph{resolvent matrix}
\[G(z)=(H-z)^{-1}, \quad \Im z>0.\]
Of particular importance is the normalized trace of $G(z)$, because it equals the Stieltjes transform of the empirical eigenvalue distribution:
\[m_N(z)=\frac{1}{N}\mathrm{tr}~G(z)=\frac{1}{N}\sum_{j=1}^N\frac{1}{\lambda_j-z}.\]
The semi-circle law is equivalent to 
\[m_N(z)=m(z)+o(1),\]
where
\[m(z)=\int \frac{1}{x-z}\rho_{\mathrm{sc}}(x)\,\mathrm{d}x = \frac{-z + \sqrt{z^2-4}}{2}.\]
The version of the local semi-circle law we state here is taken from \cite[Theorem 2.6]{KBG}, which is a detailed, pedagogical treatment of the semi-circle law and its applications.  To state it we introduce the following notion of overwhelming probability.
\bed
We say that an event or family of events $\{ \A_i \}_{i \in \I}$ hold with overwhelming probability if for all $D>0$ there is an $N_D$,  so that $\sup_{i \in \I} \pp [ A_i^c] \le N^{-D}$ for $N \ge N_D$.
\eed
\begin{thm}[Local semi-circle law]\label{thm: sc-law}
Define
\[\mathbf{S}=\{E+i\eta: |E|\le 10, 0<\eta\le 10\}.\]
Then, for each $\epsilon>0$ and each $D>0$ large, and all $N$ sufficiently large, we have
\[\max_{i,j}|G_{ij}(z)-\delta_{ij}m(z)|\le \sqrt{\frac{\Im m(z)}{N^{1-\epsilon}\eta}}+\frac{1}{N^{1-\epsilon}\eta},\]
and
\[|m_N(z)-m(z)|\le \frac{1}{N^{1-\epsilon}\eta}\]
uniformly in $z\in \mathbf{S}$ with overwhelming probability.
\end{thm}
One consequence of the semi-circle is that the eigenvalues $\lambda_i$ are close to the classical location \eqref{eqn: classical-loc}.
\begin{thm}[Eigenvalue rigidity] For each $\epsilon>0$ and $D>0$, we have
\[|\lambda_i-\gamma_i|\le N^{-2/3+\epsilon}\min\{i, (N+1-i)\}^{-1/3}\]
uniformly in $i$ with overwhelming probability.
\end{thm}

\section{Homogenization for DBM} \label{sec:homog}

\subsection{Wigner matrices} \label{sec:wighomog}

In this section we present the homogenization result of \cite{fixed1} for Wigner matrices.  We will need first the following definition.
\bed \label{def:gde}
We say that $X$ is a Wigner matrix with Gaussian component of size $t_0$ if,
\begin{align}
X = \e^{-t_0/2} X' + \sqrt{1-\e^{-t_0} } W
\end{align}
where $X'$ is a Wigner matrix and $W$ is an independent GOE.
\eed
Given a Wigner matrix $X$ (which we will eventually take to have a Gaussian component) and a GOE matrix $W'$ independent of $X$, we define the following system of coupled SDEs.  First we define,
\begin{align}
\d x_i (t) = \sqrt{ \frac{2}{N}} d B_i + \frac{1}{N} \sum_{j \neq i }\frac{1}{ x_i - x_j } \d t - \frac{x_i}{2} \d t, \qquad x_i (0) = \lambda_i (X)
\end{align}
where the $B_i$ are standard Brownian motions.  Using the same Brownian motion terms, we define
\begin{align}
\d y_i (t) = \sqrt{ \frac{2}{N}} d B_i + \frac{1}{N} \sum_{ j \neq i } \frac{1}{y_i -y_j } \d t - \frac{y_i}{2} \d t, \qquad y_i (0) = \lambda_i (W').
\end{align}
It is well-known that for each time $t$, the vector $x_i(t)$ is distributed as the eigenvalues of the matrix $\e^{-t/2}X + \sqrt{1-\e^{-t}} W''$ where $W''$ is an independent GOE matrix (see, e.g., the classic calcuation of Dyson \cite{D} or Section 4.3 of \cite{AGZ} for a pedagogical treatment).  Consequently, the vector $y_i (t)$ is distributed as the eigenvalues of a GOE matrix for every $t$.

The following homogenization result is Theorem 3.2 of \cite{fixed1}.
\bet \label{thm:wighomog}  Let $\kappa >0$ and fix an index $i$ with $\kappa N \le i \le (1-\kappa)N$.  
There are constants $\tau_0 <1/4$ and $\delta_1, \delta_2$ so that the following.  Suppose that $X$ is a Wigner matrix with Gaussian component of size $t_0 = N^{-\tau_0}$.  Let $W$ be a GOE matrix independent of $X$ and consider the coupled system of SDEs defined above.  Then, with probability at least $1-N^{-\delta_1}$ and $t_1 = t_0/2$,
\begin{align}
x_i (t_1) - y_i (t_1) = \frac{1}{N} \sum_{j} p_{t_1} (\gamma_i, \gamma_j ) (x_j (0)  - y_j (0)) + \O \left(\frac{1}{N^{1+\delta_2}} \right).
\end{align}
The function $p_{t_1} (x, y)$ is smooth and its properties are given below.
\eet

\remark We could have also used instead the homogenization result of \cite{fixed2}, which would have yielded the above estimate with much higher probability.  However, the result of \cite{fixed2} is presented in a general setting, and it would take some exposition to specialize it to the simpler Wigner case.  Moreover, we will not need the stronger result proved there.

The function $p_{s} (x, y)$ is defined on $[-2, 2]^2$ is defined explicitly in (3.22) of \cite{fixed1}.  It is smooth and obeys the estimates \cite{fixed1},
\beq \label{eqn:psbds}
0 \le p_s( \gamma_i , \gamma_j ) \le \frac{Cs}{s^2 + (\gamma_i - \gamma_j)^2}, \qquad \left| \partial_x p_s (\gamma_i, x) \right| \le \frac{C}{s^2 + (\gamma_i - x)^2}
\eeq
In preparation for the application of the mesoscopic central limit theorem, we show how to rewrite the sum on the RHS as a linear statistic.  Choose an $\eps_1 >0$ so that $0 < \eps_1 < \tau_0/10$.   Using rigidity, and the first estimate  of \eqref{eqn:psbds} we see that, 
\beq
\frac{1}{N} \sum_{j} p_{t_1} (\gamma_i, \gamma_j ) (x_j (0)  - \gamma_j )  = \frac{1}{N} \sum_j p_{t_1} ( \gamma_i, \gamma_j ) \chi \left( \frac{ \gamma_j - \gamma_i }{t_1 N^{\eps_1}} \right) ( x_j (0) - \gamma_j )+ \O \left( \frac{1}{N^{1+\eps_1/5}} \right)
\eeq
with overwhelming probability.  Here, $\chi$ is a smooth indicator function  that is $1$ for $|y| \le 1$ and $0$ for $|y| \ge 2$.  An explicit calculation using (3.22) of \cite{fixed1} yields,
\beq
p_{t_1} (x, y) =  \frac{1}{ \rhosc ( \gamma_i ) } \frac{t_1  \rhosc ( \gamma_i )}{ (x-y)^2 + ( t_1 \rhosc ( \gamma_i ))^2} \left( 1 + \O (N^{-\eps_1/5} \right), \qquad |x-\gamma_i |+ |y-\gamma_i | \le 2 t_1N^{ \eps_1}.
\eeq
Note that we left uncancelled the factor $\rhosc ( \gamma_i )$ to show the connection of $p_{t_1} (x, y)$ to the Poisson kernel $\frac{\eta}{(x-y)^2 + \eta^2}$. 
Defining,
\beq\label{eqn: Gexp}
\G(y) = \int_{-2}^{y} \chi \left( \frac{x - \gamma_i }{t_1 N^{\eps_1}} \right)    \frac{1}{ \rhosc ( \gamma_i ) } \frac{t_1  \rhosc ( \gamma_i )}{ (x-\gamma_i)^2 + ( t_1 \rhosc ( \gamma_i ))^2} d x,
\eeq
we see that, by a similar argument as Section 4 of \cite{fixed1} (near (4.57-4.58)) that on the event that the conclusion of Theorem \ref{thm:wighomog} holds, that
\beq \label{eqn:hzeta}
x_i (t_1) - y_i (t_1)  = \frac{1}{N} \zeta_x - \frac{1}{N}\zeta_y + \O ( N^{-1-c})
\eeq
for some $c>0$, 
where,
\beq
\frac{1}{N} \zeta_x = \frac{1}{N} \sum_{j} \G ( x_j )  - \G (\gamma_j ).
\eeq
For later use, we note that for $|x-\gamma_i | > t_1 N^{\eps_1}$,
\beq \label{eqn: approximate-indicator}
\G(x) = \frac{1}{ \rhosc ( \gamma_i ) }\1_{ \{ x \ge \gamma_i \}} + o (1)
\eeq
as well as the bounds,
\beq \label{eqn:Gbds}
\G(x) \le C, \qquad ||\G' (x) ||_1 \le C, \qquad ||\G''(x) ||_1 \le C \frac{N^{\eps_1}}{t_1} \ll N
\eeq

\subsection{$\beta$-ensembles}

In this section we discuss the analog of Theorem \ref{thm:wighomog} for $\beta$-ensembles.  Let us consider the following set-up.  Fix an index $i_0$ with $\kappa N \le i_0 \le (1- \kappa)N$ with $\kappa >0$, and two $\beta$-ensembles with potentials $V$ and $U$.  Denote their classical eigenvalue locations by $\gV_i$ and $\gU_i$, and equilibrium measure by $\rhoV$ and $\rhoU$.  Assume that,
\begin{align}
\gV_{i_0} = \gU_{i_0} = 0, \qquad \rhoV (0) = \rhoU(0) = \rhosc (0).
\end{align}
Let, $x_i(t)$ and $y_i(t)$ be the following coupled process,
\begin{align}
\d x_i (t) = \sqrt{ \frac{2}{N \beta}} d B_i + \frac{1}{N} \sum_{j \neq i } \frac{1}{x_i - x_j } \d t - \frac{V' (x_i )}{2} \d t
\end{align}
and
\begin{align}
\d y_i (t) = \sqrt{ \frac{2}{N \beta}} d B_i + \frac{1}{N} \sum_{j \neq i } \frac{1}{y_i - y_j } \d t - \frac{U' (y_i )}{2} \d t
\end{align}
where the $B_i$ are standard Brownian motions, and the initial data are independent $\beta$-ensembles for potentials $V$ and $U$, respectively.  Note that the distribution of each of the sets of particles is invariant under $t$ (of course the joint distribution of all $2N$ particles together is not).

The following theorem is consequence of Section 8, specifically Theorem 8.15, of \cite{fixed2}.
\bet \label{thm:betahomog}
Let $t_1 = N^{\om_1}/N$ with $0 < \om_1 < 1/10$.  Then there are $\eps, \delta_1, \delta_2 >0$ so that the following estimate holds with probability at least $1- N^{-\delta_1}$.  We have,
\begin{align} \label{eqn:betahomog}
x_{i_0} (t_1) - y_{i_0} (t_1) = \frac{1}{N} \sum_{|i_0 - j | \le N t_1 N^{\eps}} \hatp_{t_1} ( \gf_{i_0} , \gf_j ) (x_j (0) - y_j (0) ) + \O ( N^{-1-\delta_2} ).
\end{align}
\eet
In the above theorem, the function $\hatp$, as the theorem is written above, is not exactly the same as $p$ appearing earlier, however it also obeys the estimates \eqref{eqn:psbds}.  Above $\gf_j$ are the `flattened' classical eigenvalue locations,
\beq
\gf_i = \frac{i}{N \rhosc (0)}.
\eeq
Note that for $|j-i_0| \le \sqrt{N}$, the quantities $\gf_j$, $\gU_j$ and $\gV_j$ all differ by less than $C/N$.  Similar arguments as to those in the Wigner case show that the sum on the RHS of \eqref{eqn:betahomog} can be written as
\beq
 \frac{1}{N} \sum_{|i - j | \le N t_1 N^{\eps}} \hatp_{t_1} ( \gf_i , \gf_j ) (x_i (0) - y_i (0) )  = \frac{1}{N} \left( \hatzeta_x - \hatzeta_y  \right)+ \O ( N^{-1-c} )
\eeq
where
\beq
\hatzeta_x = \frac{1}{N} \sum_j  \hat{\G} ( x_i ) - \hat{\G} ( \gV_i )
\eeq
and similarly for $\hatzeta_y$.   Here $\hat{\G}$ obeys also \eqref{eqn: approximate-indicator} (with $\rhosc ( \gamma_i)$ replaced by $\rhosc (0)$) and \eqref{eqn:Gbds}.

\section{Linear statistics}
In this section, we derive the mesoscopic central limit theorem we will use to prove our main results.  Let $H$ be a real symmetric Wigner matrix, and let $a_k$ and $s_k$ be as in \eqref{eqn:cumulants}.  We will use frequently the local semicircle law discussed in Section \ref{sec:lsc}.
\subsection{Statements}\label{sec: statements}
Let $f_N (x) \in C^3$ be a sequence of test functions. We will drop the $N$-dependence from the notation and write $f = f_N$. We assume that there are $c, C>0$ so that
\beq\label{eqn: testf-assump}
\|f''\|_{\Lone}  \le N^{1-c}, \qquad \|f'\|_{\Lone} + \|f\|_{\Lone} \le C,
\eeq
where $||f||_p$ denotes the $L^p$ norm.  
We also assume that
\beq
f'(x) \neq 0 \mbox{ only if } x \in (-2 + \kappa, 2 - \kappa)
\eeq
for some $\kappa >0$.  For definiteness we assume $\kappa <\frac{1}{10}.$
 Define, for $\lambda$,
\beq
\psi ( \lambda ) = \ee[ e ( \lambda ) ] , \qquad e ( \lambda ) = \exp \left\{ \i \lambda ( \tr f(H) - \ee[ \tr f (H) ] ) \right\}.
\eeq
The function 
$\psi$ is the characteristic function of the centered linear statistic.   Let $\chi(y)$ be a smooth function that is $1$ for $|y| \le 1$ and $0$ for $|y| \ge 2$.  We use the quasi-analytic extension of $f$,
\beq
\tilf (x+ \i y) := \chi(y) ( f(x) + \i y f'(x) ).
\eeq

In this section we make the following calculation.  
\bep \label{prop:stein}  We have for any $\eps >0$,
\begin{align}\label{eqn: cf}
\del_ \lambda \psi ( \lambda ) = - \lambda \psi ( \lambda) V(f) +  \O \left(  N^{\eps} (1 + | \lambda|^5 )  N^{-1/2}(1+  \|f''\|_1^{1/2}) \right).
\end{align}
where
\begin{align}
V ( f) := - \frac{ 1 }{\pi^2} \int\int_{\cc^2} \bar{\del}_z \tilf (z) \bar{\del}_{z'} \tilf(z') \frac{1}{ z + 2 m (z) }& \bigg\{ 2 \del_{z'} \frac{ m (z') - m (z)}{z'-z} \notag\\
&+ 2 s_4 m^2 (z) m(z') m' (z')  +(a_2-1)m'(z') m (z) \bigg\} \d z \d z'
\end{align}
\eep
We can further calculate the variance and expectation.
\bel \label{lem:var}
We have,
\begin{align}\label{eqn:EV}
V(f) =& \frac{1}{2 \pi^2} \int_{-2}^2 \int_{-2}^2 \frac{ ( f(y) - f(x) )^2}{(x-y)^2} \frac{ 4 - xy}{ \sqrt{4-x^2} \sqrt{4 - y^2 } }\,\mathrm{d}x\mathrm{d}y \notag\\
&+ \frac{a_2-1}{4 \pi^2} \left( \int_{-2}^2 f(x) \frac{x}{\sqrt{4-x^2} }\,\mathrm{d}x \right)^2 \notag\\
&+ \frac{s_4}{2 \pi^2} \left( \int_{-2}^2 f(x) \frac{ 2-x^2}{\sqrt{4-x^2} } \,\mathrm{d}x \right)^2.
\end{align}
For the expectation we have for any $\eps >0$,
\begin{align}\label{eqn:Etr}
\ee[ \tr f (H) ] &= N \int_{-2}^2 f(x) \rho_{\mathrm{sc} }(x) \d x + \frac{-1}{ 2 \pi} \int_{-2}^2 \frac{f(x)}{ \sqrt{4-x^2}} \d x + \frac{ f(2) + f(-2)}{4}  \notag\\
&+\frac{1-a_2}{ 2 \pi } \int_{-2}^2 f(x) \frac{2-x^2}{ \sqrt{4-x^2}}  \d x\notag\\
&+ \frac{s_4}{2 \pi} \int_{-2}^2 f(x) \frac{ x^4-4x^2+2}{\sqrt{4-x^2}} \d x \notag\\
&+ \O\left( N^{\eps} N^{-1/2} ( \|f''\|_1^{1/2} +1 ) \right).
\end{align}
\eel
\remark The form of the variance \eqref{eqn:EV} has appeared before in the literature.  In the case that the Wigner ensemble matches the GOE to four moments, $a_2=1$ and $s_4 = 0$, and the formula matches that of Lytova and Pastur \cite{lytovapastur}.  In the general case, this formula appeared in the work of Shcherbina \cite{shcherbina}.

Integrating \eqref{eqn: cf} gives us the following theorem.
\bet \label{thm:main-meso-clt}
Let $f$ be as above and assume $V(f) \ge c$.  Then
\beq
\frac{ \tr f(H) - \ee[ \tr f(H) ] }{V(f) } 
\eeq
converges to a Gaussian with mean $0$ and variance $1$.
\eet

\subsection{Proof of Proposition \ref{prop:stein}}

We recall here the Hellfer-Sj\"ostrand formula (see \cite{davies} or Appendix C of \cite{KBG}) which states that,
we have,
\begin{align}
\tr f (H) - \ee[ \tr f (H) ] &= \frac{1}{ \pi} \int_{\cc^2} ( \del_{\bar{z}} \tilde{f} (z) ) \d z \d \bar{z} \notag \\
&= \frac{1}{2  \pi} \int_{\rr^2} \left( \i y \chi (y) f''(x) + \i ( f (x) + \i  f'(x) y ) \chi' (y) \right) N \left( m_N (z)  - \ee[ m_N (z) ] \right) \d x \d y,
\end{align}
where $\tilde f (x + \i y ) := (f (x) + \i y f' (x) ) \chi (y) $ is a quasi-analytic extension of $f$ and $\chi(y)$ is a smooth cut-off function that is $1$ for $|y| <1$ and $0$ for $|y| > 2$.  
Fix a small $\mfa > 0$ and define the domain 
\beq
\Oma := \{ (x,y ) \in \rr^2 : |y| > N^{ \mfa -1} \}.
\eeq
Using the fact that $ y \to \Im [ m_N (x + \i y) ] y $ is an increasing function and that $m_N ( \bar{z} ) = \bar{m}_N (z )$ we get,
\begin{align} 
\tr f (H) - \ee[ \tr f (H) ] &= \frac{1}{2 \pi} \int_{\Oma} \left( \i y \chi (y) f''(x) + \i ( f (x) + \i f'(x) y ) \chi' (y) \right) N \left( m_N (z)  - \ee[ m_N (z) ] \right) \d x \d y \notag\\
&  + \O\left( N^{2 \mfa-1} ||f''||_1 \right) \label{eqn:omared}
\end{align}
with overwhelming probability.  Using this, we have
\begin{align}
\psi' ( \lambda ) &= \frac{\i}{ 2 \pi} \int_{\Oma }\big( \i y \chi(y) f''(x) + \i ( f(x) + \i y f'(x) ) \chi' (y)\big) E(z) \d x \d y  + \O\left( N^{2 \mfa-1} ||f''||_1 \right)
\end{align}
where
\beq
E(z) = N \ee[ e (\lambda) (m_N (z) - \ee[ m_N (z) ]) ] = \sum_i \ee[ e ( \lambda) ( G_{ii} - \ee[G_{ii} ] )].
\eeq
The following lemma will be used repeatedly.
\bel \label{lem:sbd}
Let $H(z)$ be a holomorphic function on $\cc \backslash \rr$.  Suppose that the estimate
\beq
|H(z) | \le \frac{ K}{ |\Im [z ] |^s},
\eeq
for some $1 \le s \le 2$ whenever $|z-2| > \kappa/2$ or $|z+2| > \kappa/2$, where $\kappa$ is as in the definition of $f$. 
There is a $C>0$ so that,
\begin{align}
\left| \int_{ \Oma }\big( \i y \chi(y) f''(x) + \i ( f(x) + \i y f'(x) ) \chi' (y)\big) H (x+ \i  y ) \d x \d y \right| \le C K \log(N) \left(1+ ||f''||_1 \right)^{s-1}.
\end{align}
\eel
\proof We have only to deal with the term involving $f''(x)$.  For this, we fix a scale $\eta_1$ and integrate in parts for $|y| > \eta_1$, to find
\beq
\left| \int_{|y|> \eta_1 } y \chi (y) f''(x) H(x + \i y ) \d x \d y \right| = \left| \int_{|y| > \eta_1} \chi (y) f' (x)  y(\del_z H) (x + \i  y )  \d x \d y \right|
\eeq
where used the Cauchy-Riemann equations.  Since $H$ is holomorphic away from the real axis we have by the Cauchy integral formula,
\beq
| \del_z H (x + \i y ) | \le \frac{C K}{ | \Im [z]|^{s+1}}.
\eeq
Hence, 
\begin{align}
\left| \int_{\Oma} y \chi (y) f''(x) H \d x \d  y \right|  &\le \left| \int_{|y| > \eta_1} \chi (y) f' (x)  y(\del_z H) (x + \i  y )  \d x \d y \right| + \left| \int_{\Oma \cap \{ |y| \le \eta_1 \} } y \chi (y) f''(x) H \d x \d y \right|  \notag\\
&\le CK \log(N) \left( ||f'||_1 (\eta_1)^{1 - s} + ||f''||_1 ( \eta_1)^{2-s} \right) \notag\\
&\le C' K \log(N) \left( (\eta_1)^{1-s} + ||f''||_1 ( \eta_1)^{2-s} \right),
\end{align}
where we used the assumption that $||f'||_1 \le C$.  The result follows by taking $\eta_1 = ||f''||_1^{-1}$. \qed

\remark Note the integration by parts technique in the above proof; before the partial integration, the integral over $x$ contributes $||f''||_1$.  Integrating by parts contributes $\O ( |y|^{-1} )$ from the derivative of $H(z)$, which is smaller when $|y| > ||f''||_1^{-1}$.   The integral in the variable $y$ is always estimated by power counting. This argument of integration by parts and power counting will be used repeatedly in what follows. \qed

Let now 
\beq
\ea (\lambda) = \exp\left[ \frac{\i \lambda}{\pi} \int_{\Oma} \left( \i y \chi (y) f''(x) + \i ( f (x) + \i f'(x) y ) \chi' (y) \right) N \left( m_N (z)  - \ee[ m_N (z) ] \right) \d x \d y \right],
\eeq
so that by \eqref{eqn:omared},
\beq \label{eqn:eae}
\left| e ( \lambda) - \ea ( \lambda ) \right| \le C | \lambda | N^{2 \mfa -1} ||f''||_1.
\eeq
Let $\Ea(z)$ be the same as $E(z)$ but with $\ea (\lambda)$ in place of $e ( \lambda)$.  Applying Lemma \ref{lem:sbd} with $H = E - \Ea = \O ( |\lambda| N^{\eps+ 2 \mfa -1} ||f''||_1 \eta^{-1})$,  we see that
\beq \label{eqn:psidera}
\psi' ( \lambda) = \frac{ \i}{ \pi} \int_{\Oma} ( i y \chi (y) f'' (x) + \i  ( f (x) + \i y f' (x) ) \chi' (y) ) \Ea (z) \d x \d y + \O ( N^{\eps+2 \mfa -1} ||f''||_1 )
\eeq
for any $\eps >0$. 

We will use a cumulant expansion on the right-hand side of \eqref{eqn:psidera}.  Let $\del_{ab}$ denote differentiation with respect to $H_{ab}$.  By a direct calculation (similar to \eqref{eqn:eafirstder} below)  using Lemma \ref{lem:sbd} one can see that with overwhelming probability,
\beq \label{eqn:dereabd}
| \del_{ab}^k \ea ( \lambda ) | \le (1 + | \lambda | )^k  N^{\eps }
\eeq
for any $\eps >0$ and $N$ large enough.
\bel  \label{lem:cumulant}
For $z \in \Oma$ we have,
\begin{align}
(z+2 m(z) ) \Ea (z) &= \sum_{i a} \frac{1+a_2 \delta_{ia}}{N} \ee[ \del_{ia} \ea ( \lambda )  ( G_{ia} )] -  \sum_{a \neq i } \frac{ s_4}{2 N^2} \ee[ \del_{ai}^2 \ea ( \lambda) m^2 ]  \notag\\
&+\O\left( N^{\eps} ( 1 + |\lambda|^4 )(  (N \eta)^{-1/2} ( 1 + \|f''\|_{\Lone} )+ (N \eta^2)^{-1} + N^{-1/2} \eta^{-1} ) \right) \label{eqn:cumulantresult}
\end{align}
\eel
\proof By the cumulant expansion (see, e.g., Lemma 3.2 of \cite{LS} and the discussion therein for use of the cumulant expansion),  we derive
\begin{align}
&z \ee[ \ea ( \lambda)  (G_{ii} - \ee[G_{ii} ] ) ]  = \sum_a \ee[ \ea ( \lambda ) ( G_{ia} H_{ai} - \ee[ G_{ia} H_{ai} ] ) ]  \notag\\
&= \frac{1+\delta_{ia} a_2}{N} \sum_{a } \ee[  \del_{ia} \ea ( \lambda )  G_{ia} ] \label{eqn:l1} \\
&- \frac{1+\delta_{ia} a_2}{N} \sum_{a } \ee[ \ea ( \lambda ) ( G_{ii} G_{aa} - \ee[ G_{ii} G_{aa} ] ) ] \label{eqn:l2} \\
&- \frac{1-\delta_{ia} }{N} \sum_{a} \ee[ \ea ( \lambda) ( G_{ia}^2 - \ee[ G_{ia}^2] )] \label{eqn:l3} \\
&+ \frac{s_3+ \delta_{ia} a_3}{2 N^{3/2}} \sum_a \ee[ \del_{ai}^2 \ea ( \lambda) G_{ai} ] + 2\ee[ \del_{ai} \ea ( \lambda)  \del_{ai} G_{ai} ] + \ee[ \ea ( \lambda) ( \del_{ai}^2 G_{ai} - \ee[ \del_{ai}^2 G_{ai} ] )] \label{eqn:l4}\\
&+ \frac{s_4+\delta_{ia} a_4}{6 N^2} \sum_a \ee[ \del_{ai}^3 \ea G_{ai} ] + \cdots + \ee[ \ea ( \lambda) ( \del_{ai}^3 G_{ai} - \ee[\del_{ai}^3 G_{ai} ])]  \label{eqn:l5} \\
&+ \O \left( N^{\eps-{3/2}} (1+ |\lambda|^4) \right),
\end{align}
any $\eps >0$. 
We begin calculating each term.  Starting with \eqref{eqn:l2} we have,
\begin{align}
\sum_{i,a} \frac{1+\delta_{ia} a_2}{N} \ee[ \ea ( \lambda ) ( G_{ii} G_{aa} - \ee[ G_{ii} G_{aa} ] ) ] &= 2 \sum_i m \ee[ \ea ( \lambda) ( G_{ii} - \ee[G_{ii} ] ) ) ] \notag\\
&+ \O ( N^{\eps} ( ( N \eta)^{-1/2} + (N \eta^2)^{-1} )).
\end{align}
This term contributes the $2m \Ea$ to the LHS of \eqref{eqn:cumulantresult}. 
For \eqref{eqn:l3} we have,
\begin{align}
\sum_{i,a} \frac{1-\delta_{ia} }{N}  \ee[ \ea ( \lambda) ( G_{ia}^2 - \ee[ G_{ia}^2] )] &= \ee[ \ea ( \lambda) \del_z (m_N - \ee[ m_N ] ) ] + \O ( N^{\eps} ( N \eta)^{-1/2}) \notag\\
&= \O \left( N^{\eps} ( ( N \eta)^{-1/2}  + (N \eta^2 )^{-1} ) \right).
\end{align}
We now start with the terms in \eqref{eqn:l4}.  For a combinatorial factor $K$,
\begin{align}
\sum_{i,a} \frac{s_3 + \delta_{ia} a_3}{N^{3/2}} \ee[ e ( \lambda) ( \del_{ai}^2 G_{ai} - \ee[ \del_{ai}^2 G_{ai} ] )] &= K m^2 \frac{1}{N^{3/2}} \sum_{i a } \ee[ e ( \lambda) (G_{ia} - \ee[ G_{ia} ] )] \notag\\
&+ \O ( N^{\eps} ( N^{-1/2} + N^{-1/2} \eta^{-1} ))
\end{align}
By the isotropic local law (Theorem 2.12 of \cite{BEKYY}),
\beq
\left| \frac{1}{N^{1/2}} \sum_i G_{ia} \right| \le \frac{N^{\eps}}{\sqrt{N \eta}}
\eeq
with overwhelming probability, and so
\beq
\frac{1}{N^{3/2}} \sum_{i a } \ee[ e ( \lambda) (G_{ia} - \ee[ G_{ia} ] )] = \O ( N^{\eps} (N \eta)^{-1/2} ).
\eeq
In preparation for the second term of \eqref{eqn:l4} we calculate, for $i \neq a$,
\begin{align} \label{eqn:eafirstder}
\del_{ia} \ea ( \lambda ) = - \frac{2 \i \lambda}{\pi} \ea ( \lambda) \int_{\Oma} \bar{\del}_{z' } \tilf ( z' ) \sum_b G_{bi}G_{ab} = -\frac{2 \i \lambda}{\pi} \ea ( \lambda) \int_{\Oma} \bar{\del}_{z'} \tilf ( z') \del_{z'} G_{ia}.
\end{align}
Since with overwhelming probability,
\beq
\sum_b G_{bi}G_{ba} = \del_{z'} G_{ai} = \O ( N^{\eps} (N \eta')^{-1/2} (\eta')^{-1} ),
\eeq
we find that (using Lemma \ref{lem:sbd})
\beq \label{eqn:firstder}
| \del_{ia} \ea ( \lambda ) | \le N^{\eps} (1 + |\lambda| ) \left( N^{-1/2} ( 1+||f''||_{1})^{1/2}   )\right).
\eeq
For the second term of \eqref{eqn:l4} we have
\begin{align}
\sum_{ia} \frac{s_3+ \delta_{ia} a_3}{N^{3/2}} \ee[ \del_{ai} \ea ( \lambda)  \del_{ai} G_{ai} ] = - \frac{s_3 }{N^{3/2}} \sum_{i \neq a} \ee[ \del_{ai} \ea ( \lambda ) G_{aa} G_{ii} ] + \O ( N^{\eps} (1+|\lambda|)(N^{-1/2} + N^{-1/2} \eta^{-1} )).
\end{align}
We see, using \eqref{eqn:firstder},
\beq
 \frac{s_3 }{N^{3/2}} \sum_{i \neq a} \ee[ \del_{ai} \ea ( \lambda ) G_{aa} G_{ii} ]  =  \frac{s_3 m^2}{N^{3/2}} \sum_{i \neq a} \ee[ \del_{ai} \ea ( \lambda )  ]  + \O \left( N^{\eps} \eta^{-1/2} N^{-1/2} (1 + \|f'\|_{\Lone}^{1/2} ) (1 + | \lambda | ) \right).
\eeq
The derivative in the above expression together with the sum over $i, a$ gives us $ N^{-3/2} \del_{z'} \sum_{i,a} G_{ia} = N^{-1/2} \del_{z'} m + \O ( N^{-3/2} (\eta')^{-2} )$ by the isotropic local law.  So, using Lemma \ref{lem:sbd},
\beq
\frac{1 }{N^{3/2}} \sum_{i \neq a} \ee[ \del_{ai} \ea ( \lambda ) ]  = \O ( N^{\eps} ( 1 + |\lambda| )N^{-1/2} ).
\eeq
For the first term of \eqref{eqn:l4} we further calculate, for $i \neq a$,
\begin{align}
\del_{ia}^2 \ea ( \lambda) &= \ea ( \lambda)  (  i \lambda)^2 \pi^{-2} \left( \int_{\Oma} \bar{\del}_{z'} \tilf ( z' ) \sum_b G_{bi} G_{ab} \d z' \d \bar{z}' \right)^2 +2 \ea ( \lambda) \frac{ \i \lambda}{\pi} \int_{\Oma} \bar{\del}_{z'} \tilf (z') \del_{z'} (G_{aa} G_{ii}  + G_{ia}^2) \notag\\
&= \ea ( \lambda) \frac{ \i \lambda }{\pi} \int_{\Oma} \bar{\del}_{z'} \tilf ( z' ) 4 m(z') m' (z' ) \d z' \d \bar{z}'+ \O \left( N^{\eps} ( 1 + | \lambda | )^2 N^{-1/2} ( \|f''\|_{\Lone}^{1/2} +1 )\right) \label{eqn:der2}
\end{align}
and so 
\begin{align}
\sum_{a, i} \frac{s_3+ \delta_{ia} a_3}{N^{3/2}}  \ee[ \del_{ai}^2 \ea ( \lambda) G_{ai} ] &= \sum_{a, i} \frac{s_3}{N^{3/2}} \ee\left[ \ea ( \lambda) \i \lambda \pi^{-1} \int_{\Oma} \bar{ \del}_{z'} \tilf (z') 4 m (z') m' (z') G_{ai} \d z' \d \bar{z}' \right] \notag\\
&+ \O \left( N^{\eps} ( 1 + | \lambda | )^2 (N \eta)^{-1/2} ( \|f''\|_{\Lone}^{1/2} +1)\right) \notag\\
&=  \O \left( N^{\eps} ( 1 + | \lambda | )^2 (N \eta)^{-1/2} ( \|f''\|_{\Lone}^{1/2} + 1)\right) 
\end{align}
where in the last line we used the isotropic law again.

Finally, for the terms \eqref{eqn:l5} we see that (using \eqref{eqn:dereabd} and the Green's function estimates of Theorem \ref{thm: sc-law} repeatedly)
\begin{align}
 &\sum_{a,i} \frac{s_4+\delta_{ia} a_4}{N^2}  \ee[ ( \del_{ai}^3 \ea)  G_{ai} ] + \cdots + \ee[ \ea ( \lambda) ( \del_{ai}^3 G_{ai} - \ee[\del_{ai}^3 G_{ai} ])] \notag\\
 = & -  \sum_{a \neq i } \frac{3 s_4}{N^2} \ee[ \del_{ai}^2 \ea ( \lambda) m^2 ] + \O\left( N^{\eps} (1 + |\lambda|^3) (N \eta)^{-1/2} \right)
\end{align}
\qed

We further calculate the terms on the right side of \eqref{eqn:cumulantresult}.

\bel 
We have,
\begin{align}
\sum_{i a} \frac{1+  \delta_{ia}}{N} \ee[ ( \del_{ia} \ea ( \lambda ) ) ( G_{ia} )]  = - \frac{2 \i \lambda}{\pi} \int_{\Oma} \bar{\del}_{z'} \tilf (z') \del_{z'} \frac{1}{N}\ee [ e ( \lambda) \tr G (z') G(z)] \d z' \d \bar{z}'.
\end{align}
and
\begin{align} \label{eqn:deliie}
\frac{1}{N} \sum_i \ee[ \del_{ii} \ea ( \lambda) G_{ii} ] = -  \frac{ \i \lambda \ee[ \ea ( \lambda) ]}{\pi} \int_{\Oma} \bar{\del}_{z'} \tilf (z') m' (z') m (z) \d z' \d \bar{z}' + \O \left( N^{\eps} (1 + | \lambda ) ) (N \eta)^{-1/2} \right)
\end{align}
and
\begin{align}
&\sum_{a \neq i } \frac{ s_4}{2 N^2} \ee[ \del_{ai}^2 \ea ( \lambda) m(z)^2 ]  \notag\\
=& 2 s_4 \ee[ \ea ( \lambda)] \i \lambda \int_{\Oma} \bar{\del}_{z'} \tilf ( z' )  m(z') m' (z' ) m (z)^2 \d z' \d \bar{z}' + \O \left( N^{\eps} ( 1 + | \lambda | )^2 N^{-1/2} ( \|f''\|_1^{1/2} +1 )\right)
\end{align}
\eel
\proof We see that,
\begin{align}
\sum_{i a} \frac{1+\delta_{ia}}{N} \ee[ ( \del_{ia} e ( \lambda ) ) ( G_{ia} )] &= -\frac{ 2 \i \lambda}{\pi N } \int \bar{\del}_{z'} \tilf (z') \ee \sum_{iab} e ( \lambda) G_{ab} (z') G_{bi} (z') G_{ia} (z) \d z' \d \bar{z}' \notag\\
&= -\frac{ 2 \i \lambda}{\pi} \int \bar{\del}_{z'} \tilf (z') \del_{z'} \frac{1}{N}\ee \left[ e ( \lambda) \tr G (z') G(z)\right] \d z' \d \bar{z}'.
\end{align}
The second estimate is similar, and the last estimate follows from \eqref{eqn:der2}. 
\qed

Putting together the last two lemmas we obtain,
\begin{align}
\Ea ( z) &= \frac{1}{ z + 2 m (z) } \frac{ \i \lambda}{\pi} \int_{\Oma}  \d z' \d \bar{z}' \bar{\del}_{z'} \tilf (z')  \notag\\
\times  & \bigg\{ \del_{z'} \frac{-2}{N} \ee[ \ea ( \lambda)  \tr G(z') G(z) ] - \ee[\ea ( \lambda)] (2 s_4m(z)^2 m(z') \del_{z'} m(z')-(a_2-1) m'(z') m (z)) \bigg\}   \notag\\
&+ \O\left( N^{\eps} ( 1 + |\lambda|^4 )(  (N \eta)^{-1/2} ( 1 + \|f''\|_{L^1(\mathbb{R})}^{1/2} )+ (N \eta^2)^{-1} + N^{-1/2} \eta^{-1} ) \right)
\end{align}

We now complete the proof of Proposition \ref{prop:stein}.   We see that with overwhelming probability, if $z$ and $z'$ are on different half-planes,
\begin{align}
\del_{z'} \frac{1}{N} \tr G(z') G(z) &= \del_{z'} \frac{1}{N} \tr \frac{ G(z') - G(z) }{z'-z} \notag\\
&= \del_{z'} \frac{ m(z') - m(z)}{z'-z} + \O \left( N^{\eps} ( \eta')^{-1} ( (N \eta')^{-1} + (N \eta)^{-1} )(\eta + \eta')^{-1} \right)
\end{align}
whereas if they are on the same half-plane, divide into the cases $\eta < 2 \eta'$ and $ \eta > 2 \eta'$.  In the former, write $m(z) - m(z') = \int m' (s) \d s$ to get
\begin{align}
\del_{z'} \frac{1}{N} \tr G(z') G(z) &= \del_{z'} \frac{ m(z') - m(z)}{z'-z}  + \O \left( N^{\eps} ( \eta')^{-1} ( ( N \eta^2 )^{-1} + ( N ( \eta')^2 )^{-1} ) \right) \notag\\
&= \del_{z'} \frac{ m(z') - m(z)}{z'-z}  + \O \left( N^{\eps} ( \eta')^{-1} ( N \eta^2 )^{-1} \right).
\end{align}
Whereas if $\eta > 2 \eta'$ directly estimate $m_N (z') - m_N (z) = m(z') - m (z) + \O ( N^{\eps} ( N \eta')^{-1} ))$ to get,
\begin{align}
\del_{z'} \frac{1}{N} \tr G(z') G(z) &= \del_{z'} \frac{ m(z') - m(z)}{z'-z}   + \O \left( N^{\eps}( ( N \eta^2 \eta' )^{-1} + (N (\eta')^2 \eta )^{-1} )\right).
\end{align}
Integrating all of this (i.e., using the integration by parts and power counting technique of Lemma \ref{lem:sbd}) we see that for any $\eps >0$,
\begin{align}
&\frac{ \d }{ \d \lambda} \psi ( \lambda) \\
 = & \lambda \ee [ \ea ( \lambda ) ]  \bigg[ \frac{ 1 }{\pi^2} \int_{\Oma^2} \bar{\del}_z \tilf (z) \bar{\del}_{z'} \tilf(z') \frac{1}{ z + 2 m (z) } \notag\\
 \times &\left(2 \del_{z'} \frac{ m (z') - m (z)}{z'-z} + 2 s_4 m^2 (z) m(z') m' (z')  +(a_2-1)m'(z') m (z)\right) \d z \d \bar{z} \d z' \d \bar{z}' \bigg] \label{eqn:fs1} \\
&+ \O \left(  N^{\eps} (1 + | \lambda|^4 ) \left[  N^{-1/2} \|f''\|_1  + N^{ 2 \mfa-1} ||f''||_1 \right] \right).
\end{align}
In order to complete the proof, we need only estimate the error in first replacing $\ee[ \ea ( \lambda ) ]$ by $\ee[ e ( \lambda)]$ and then restoring to the integral  the region  $\cc^2 \backslash \Oma^2$.  First we estimate the size of the integral appearing above, and show that it is at most $\O (\log(N)^2)$.  We can then  replace $\ee[\ea]$ by $\ee[e ( \lambda ) ] = \psi ( \lambda)$ using \eqref{eqn:eae}.  The terms which need to be estimated involve $f''(x)$; the contribution from the other terms are bounded, due to the  assumptions \eqref{eqn: testf-assump} and the appearance of $\chi'(y)$ in these terms which is non-zero only for $y$ of order $1$.  We turn to the terms with $f''(x)$.  The second two terms in the line \eqref{eqn:fs1} are bounded functions, so we can integrate by parts in $x$ to estimate this contribution by $\O (1)$.

When $z$ and $z'$ are in the same half-spaces, $(m(z') - m(z) ) / (z' - z)$ is bounded, and so the derivative is bounded by $C/ \Im [ z']$.  Using then (the proof of) Lemma \ref{lem:sbd} we can estimate this contribution by $\O ( \log (N)$. When $z$ and $z'$ are in separate half-spaces we use instead the estimate 
\beq \label{eqn:final1}
\left| \del_{z'} \frac{ m (z) - m(z')}{ z - z'} \right| \le \frac{C}{ (  \Im [z'] )^2+ (\Im [ z] )^2}.
\eeq
The largest contribution is from when both $\Im[z']$ and  $\Im [z] > ||f''||_1^{-1}$, where we integrate by parts in both $\Re[z]$ and $\Re[z']$ and find a contribution of order $\O ( \log(N)^2)$.  The other regions are at most $\O ( \log(N))$.

Finally, we argue that we can restore to the integral the region $\cc^2 \backslash \Oma^2$ at an error of at most $\O (  \log(N) N^{2 \mfa -1} ||f''||_1)$.  The region $( \cc \backslash \Oma)^2$ is easy to control as the second line of \eqref{eqn:fs1} is bounded by $C/ ( \Im [z] \Im [z'])$, and gives a contribution of $ \O ( N^{2 \mfa - 2 } ||f''||_1^2)$.  The ``cross terms'' $\Oma \times ( \cc \backslash \Oma )$ (and vice versa) must be handled by integration by parts again in the region where $\Im [z]$ or $\Im [z'] > ||f''||_1^{-1}$ as above, and the cases of same and different half-spaces are treated similarly, e.g., using the estimate \eqref{eqn:final1}.   The cross terms are found to contribute $\O ( N^{\mfa-1} ||f''||_1 \log(N))$.  The claimed estimate of Proposition \ref{prop:stein} follows after taking $\mfa = \eps$.
\qed

\subsection{Proof of Lemma \ref{lem:var}}
We apply Green's theorem to each term, which states
\beq
\int_{\Omega} \bar{\del}_z F (z) \d z \d \bar{z} = \frac{-\i}{2} \int_{\del \Omega} F(z) \d z
\eeq
 We write $z' = y \pm \i 0$ and $z = x \pm \i 0$.  First, we integrate by parts to find
\begin{align}
 &\frac{ -1 }{\pi^2} \int_{\cc^2} \bar{\del}_z \tilf (z) \bar{\del}_{z'} \tilf(z') \frac{1}{ z + 2 m (z) }\left(2 \del_{z'} \frac{ m (z') - m (z)}{z'-z} \right)\d z \d \bar{z} \d z' \d \bar{z}' \notag\\
 = & \frac{1}{ \pi^2}\int_{\cc^2} \bar{\del}_z \tilf (z) \bar{\del}_{z'} \del_{\Re[z']} \tilf(z') \frac{1}{ z + 2 m (z) }\left(2 \frac{ m (z') - m (z)}{z'-z} \right)\d z \d \bar{z} \d z' \d \bar{z}' .
 \end{align}
Note that $x \pm i 0 + 2 m ( x \pm \i 0 ) = \pm \i \sqrt{4-x^2} \1_{\{ |x| \le 2 \} }$.  Hence,
\begin{align}
&\frac{2}{ x + 2 m ( x \pm \i 0 ) } \left( \frac{ m (y + \i 0) - m ( x \pm \i 0)}{ y-x} - \frac{ m (y -\i 0) - m ( x \pm \i 0)}{ y-x} \right) \notag\\
= & \pm 2  \frac{ \sqrt{ 4 - y^2}}{\sqrt{ 4 - x^2} (y-x) } .
\end{align}
Hence, we find (note that some of the integrals over the real axes require intepretation as principal values)
\begin{align}
&  \frac{1}{ \pi^2}\int_{\cc^2} \bar{\del}_z \tilf (z) \bar{\del}_{z'} \del_{\Re[z']} \tilf(z') \frac{1}{ z + 2 m (z) }\left(2 \frac{ m (z') - m (z)}{z'-z} \right) \d z \d \bar{z} \d z' \d \bar{z}' \notag\\
= &\frac{-1}{\pi^2} \int_{-2}^2 \int_{-2}^2 f (x) f' (y) \frac{ \sqrt{ 4 - y^2}}{\sqrt{ 4 - x^2}(y-x)} \d x \d y \notag\\
= &  \frac{1}{ \pi^2}\int_{-2}^2 f'(y) \sqrt{4-y^2} \int_{-2}^2 \frac{f(y) - f(x)}{\sqrt{4-x^2}(y-x) } \d x \d y \notag\\
= & \frac{-1}{ 2  \pi^2} \int_{-2}^2 \int_{-2}^2 \frac{ ( f(y) - f(x) )^2}{(y-x)^2} \frac{ 4 - xy}{ \sqrt{4-x^2} \sqrt{4 - y^2 } } \d x \d y.
\end{align}
where the last line follows by integration by parts in $y$ and the second last line from the identity
\beq
\int_{-2}^2 \frac{1}{(x-y) \sqrt{4-x^2}} \d x = 0.
\eeq
For the other terms, note that $(z+2m)m = m^2-1$ and since $m' = m^2/ (1-m^2)$,
\beq
\frac{m}{z+2m} = -m'.
\eeq
Therefore,
\begin{align}
\frac{-1}{\pi^2} \int_{\cc^2} \bar{\del}_z \tilf (z) \bar{\del}_{z'} \tilf (z') \frac{ m' (z') m(z)}{z + 2m (z) } \d z \d \bar{z} \d z' \d \bar{z}'&= \frac{1}{\pi^2} \int_{\cc^2}  \bar{\del}_z \tilf (z) \bar{\del}_{z'} \tilf (z') m'(z') m' (z)\d z \d \bar{z} \d z' \d \bar{z}' \notag\\
&= \frac{1}{ 4 \pi^2} \left( \int_{-2}^2 f(x) \frac{x}{\sqrt{4-x^2} } \d x \right)^2.
\end{align}
For the final term we have,
\begin{align}
&\frac{-1}{\pi^2} \int_{\cc^2}  \bar{\del}_z \tilf (z) \bar{\del}_{z'} \tilf (z') \frac{m^2 (z) m' (z) m (z) }{z + 2 m (z) } \d z \d \bar{z} \d z' \d \bar{z}' \notag\\
=& \frac{1}{\pi^2} \int_{\cc^2} \bar{\del}_z \tilf (z) \bar{\del}_{z'} \tilf (z') m' (z') m (z') m' (z) m (z) \d z \d \bar{z} \d z' \d \bar{z}' = \frac{1}{4  \pi^2} \left( \int_{-2}^2 f(x) \frac{2-x^2}{ \sqrt{4-x^2}} \d x \right)^2.
\end{align}
This completes the calculation of $V(f)$.  
To calculate the expectation, we use the cumulant expanson on $\ee[ G_{ii}]$.   Let $\Oma$ and $\mfa >0$ be as earlier.  For $z \in \Oma$ we find for any $\eps >0$,
\begin{align}
1+ z\ee[G_{ii}] &=- \frac{1}{N} \sum_{ a \neq i } \ee[ G_{ii} G_{aa} + G_{ia}^2] - \frac{1+a_2}{N} \ee[ m^2] \notag\\
&+\frac{s_3 m^2}{2 N^{3/2}} \sum_{a \neq i } \ee[ G_{ia} ] -  \frac{s_4}{N} m^4 \notag\\
&+ \O\left( N^{\eps} N^{-1} ( N \eta^2)^{-1/2}  \right).
\end{align}
The term with $s_3$ can be absorbed into the error term using the isotropic local law.  We therefore find,
\begin{align}
\sum_i (z + 2 m ) \ee[ G_{ii} - m ] &= -\del_z m  - (a_2-1) m^2 - s_4 m^4 \notag\\
&+ \O \left( N^{\eps}(  (N \eta^2)^{-1} + (N \eta^2)^{-1/2} ) \right).
\end{align}
Integrating, we find
\begin{align}
\ee[f(H)] - N \int f  (x) \d \rho_{\mathrm{sc}}  (x) \d x &= \frac{1}{\pi}\int \bar{ \del}_z \tilf (z) \frac{ (1-a_2) m^2  - s_4 m^4 - \del_z m }{z + 2 m } \d z \d \bar{z}' \notag\\
&+ \O \left( N^{\eps}  N^{-1/2} ( \| f''\|_{1} + 1 )  + N^{ \mfa + \eps-1} ||f''||_1 \right).
\end{align}
Note that we had to remove and then return the domain $\cc\backslash \Oma$ as in the proof of Proposition \ref{prop:stein}; this is similar to what was done above.  We take $\mfa = \eps$. 
We calculate using Green's theorem as above,
\begin{align}
\frac{1}{\pi} \int_\cc \bar{ \del}_z \tilf (z) \frac{m^2}{z + 2 m } \d z \d \bar{z}' &= -\frac{1}{\pi} \int_\cc \bar{ \del}_z \tilf (z) m' (z) m (z) \d z \d \bar{z}' \notag\\
&=  \frac{1}{ 2 \pi } \int_{-2}^2 f (x) \frac{2-x^2}{ \sqrt{4-x^2}} \d x,
\end{align}
and
\begin{align}
\frac{-1}{\pi} \int_\cc \bar{\del}_{z} \tilf (z) \frac{m^4}{z+2m} \d z \d \bar{z}'&= \frac{1}{\pi} \int_\cc \bar{\del}_z \tilf (z) m^3 (z) m' (z) \d z \d \bar{z}' \notag\\
&= \frac{1}{2 \pi} \int_{-2}^2 f(x) \frac{x^4-4x^2+2}{ \sqrt{4 - x^2}} \d x.
\end{align}
Finally,
\begin{align}
\frac{-1}{\pi} \int_\cc \bar{\del}_z \tilf (z) \frac{ m'(z)}{z + 2m } \d z \d \bar{z}' &= \frac{-1}{ 2 \pi} \int_{-2}^2 \frac{f(x)}{ \sqrt{4-x^2}}\d x+ \frac{ f(2) + f(-2)}{4} .
\end{align}
This completes the proof.
\qed 

\section{Derivation of the main results}

\subsection{Reduction to Gaussian divisible ensembles}

We will only provide full details for the proofs of Theorems \ref{thm1} and \ref{thm:evalues}  for those  Wigner matrices with Gaussian components (see Definition \ref{def:gde}).  The reduction from general Wigner matrices to those with a Gaussian component is well-known in the random matrix literature, and so we will be brief in our discussion.  In our setting, this reduction will be a consequence of the four moment method of Tao and Vu \cite{TV}.  For the specific case of single eigenvalue fluctuations, similar arguments were made in \cite{TV} and \cite{orourke}.  First, we have
\bep
Let $H$ and $W$ be two Wigner matrices so that
\beq \label{eqn:approxmatch}
\left| \ee[ H_{ij}^s] - \ee[W_{ij}^s ] \right| \le N^{-2-c}, \qquad 1 \le s \le 4
\eeq
for some $c>0$.  Then if the results of Theorems \ref{thm1} or \ref{thm:evalues} hold for $W$, they hold for $H$.
\eep
This is essentially a consequence of Theorem 15 of \cite{TV}, specialized to the real symmetric setting.  The only extension is that the original Theorem asks for equality of the first four moments rather then the fact that they are only approximately equal.  This extension is well-known, see, e.g., \cite{ERSTVY}.

The additional required input is that given a Wigner matrix, one can find a matching Gaussian divisible ensemble.  The following is a consequence of, e.g.,  Lemma 3.4 of \cite{EYY}.

\bel
Let $H$ be a Wigner matrix, and let $\tau_0 >0$.  Then there is a Wigner matrix $W$ with Gaussian component of size $N^{-\tau_0}$ so that  \eqref{eqn:approxmatch} holds with some $c>0$.
\eel

The consequence of these two results is that if there is a $\tau_0 >0$ so that we can prove our main results for Wigner matrices with Gaussian components of size $N^{-\tau_0}$, then our results extend to all Wigner matrices.  The $\tau_0$ we take is the one so that Theorem \ref{thm:wighomog} holds.

\subsection{Proof of Theorem \ref{thm1}}

 Let $F$ be a smooth function, which we moreover take to be of compact support.  Let $\bfx (t)$ and $\bfy (t)$ be as in Section \ref{sec:wighomog}, so that the initial data $\bfx (0)$ are the eigenvalues of a Wigner matrix with Gaussian component of size $N^{-\tau_0}$.   We will need to introduce a third process $\bfz(t)$ coupled to the same Brownian motions as $\bfx(t)$ and $\bfy(t)$, with initial data an independent GOE matrix. 
We let
\[ Z_i(\textbf{x})=N\frac{ x_i(H_t) - \gamma_i}{\sqrt{ \frac{\log N}{(1-\gamma_i^2/4)}}},\]
and $Z_i(\mathbf{y})$ be the corresponding quantity for a GOE matrix. 
It suffices to show
\[\mathbb{E}[F(Z_i(\bfx (t)))]=\mathbb{E}[F(Z_i( \bfz (t)))]+o(1).\]  
where the time parameter $t =t_1$ as in Theorem \ref{thm:wighomog}.  By \eqref{eqn:hzeta} we have
\[\mathbb{E}[F(Z_i(\mathbf{x} (t)))]=\mathbb{E}[F\big(Z_i(\mathbf{y} (t))+\alpha_N(\zeta_{\textbf{x}}-\zeta_{\textbf{y}})\big)]+\O(N^{-c_1}),\]
where
\[\alpha_N=\frac{2\pi}{\sqrt{\log N}} \times \sqrt{1- \gamma_i^2/4}.\]
Denoting the characteristic function of $\alpha_N \zeta_x$ by $\psi_x ( \lambda)$, we have,
\[ \mathbb{E}[F\big(Z_i( \bfy (t) )+\alpha_N(\zeta_x-\zeta_y))]=\int \psi(\lambda)\widehat{F}(\lambda)\mathbb{E}[e^{i\lambda(Z_i( \bfy (t))-\alpha_N \zeta_y)}]\,\mathrm{d}\lambda.\]
Repeating the same argument for $\bfz (t)$ we see that
\beq
\ee[ F ( Z_i (\bfx (t) ) ) ] - \ee [F ( Z_i (\bfz (t) ) )  ] = \int \left( \psi_x ( \lambda) - \psi_z ( \lambda )  \right) \hat{F} ( \lambda ) \mathbb{E}[e^{i\lambda(Z_i( \bfy (t) )-\alpha_N \zeta_y)}]\,\mathrm{d}\lambda + o (1).
\eeq
Denote by $V_x (\G)$ and $V_z (\G)$ the functionals appearing in Lemma \ref{lem:var} for matrices $X$ and $Z$.  Then we have
\beq
\left| V_x (\G) - V_z (\G) \right| \le C.
\eeq
From this and Proposition \ref{prop:stein}, we obtain that
\[\psi_x(\lambda)= e^{-\alpha \lambda^2}+o(1),\quad |\lambda|\le \log(N)^{1/4},\]
where $\alpha$ is independent of the choice of $x$ or $z$. 
Since $\hat{F} $ is a Schwartz function, this yields the claim. \qed

\subsection{Proof of Corollary \ref{cor:counting}}
The argument leading from Theorem \ref{thm1} to the Corollary \ref{cor:counting} was detailed in \cite[Theorem 1.1]{gustavsson}. It suffices to notice that
\[\mathbb{P}\left(\frac{\num(E)-N \int_{-2}^E \rho_{\mathrm{sc}}(u) \d u}{\sqrt{\frac{1}{2\pi^2}\log N}}\le x\right)=\mathbb{P}\big(\lambda_{i_N}\le E\big)=\mathbb{P}\left(N\frac{\lambda_{i_N}-\gamma_{i_N}}{\sqrt{\log N}}\le x_N\right),\]
where $i_N$ is the integer part of
\[N\int^E_{-2} \rho_{\mathrm{sc}}(y)\,\mathrm{d}y-x\sqrt{\frac{1}{2\pi^2}\log N},\]
and $x_N$ is defined by the equality.  Note that,
\begin{align*}
  \frac{i_N}{N}&\rightarrow \int^E_{-2}\rho_{\mathrm{sc}}(y)\,\mathrm{d}y.\\
  x_N&\rightarrow x.
\end{align*}
The result now follows directly from the convergence in distribution of $\lambda_{i_N}$.\qed

\subsection{Proof of Theorem \ref{thm: beta}}
Recall the following result for mesoscopic linear statistics, \cite[Theorem 6.18]{fixed2}.
\begin{thm}\label{thm: beta-clt}
Suppose $f$ is a $C^2$ function with support in $[A,B]$, where $\mathrm{supp}\rho_V=[-2,2]$, satisfying the assumptions in Section \ref{sec: statements}. Let $\lambda_i$, $1\le i\le N$ be sampled from a $\beta$-ensemble with potential $V$. Then, there exists $\epsilon>0$ such that
\begin{equation}
\mathbb{E}[e^{i\lambda(\sum_{i=1}^N f(\lambda_i)-N\int f(x)\rho_V)}]=\exp\left(-\frac{\lambda^2}{2}V(f)+i\delta(f)\right)+\O(N^{-1+10\epsilon})\|f''\|_{L^1(\mathbb{R})}+\O(N^{-\epsilon}),
\end{equation}
where
\begin{align*}
V(f)&=\frac{1}{2\beta\pi^2}\int_{A}^B\int_A^B \left(\frac{f(x)-f(y)}{x-y}\right)^2\frac{-AB-xy+\frac{1}{2}(A+B)(x+y)}{\sqrt{(x-A)(B-x)}\sqrt{(y-A)(B-y)}}\,\mathrm{d}x\mathrm{d}y,\\
\delta(f)&=\frac{1}{2\pi^2}\left(\frac{2}{\beta}-1\right)\int_A^B f(x)  \d \nu_V (x)
\end{align*}
where $\nu_V$ is a signed measure.
\end{thm}
\remark The measure $\nu_V$ is characterized explicitly in the book of Pastur and Shcherbina \cite{ps} (see Theorem 11.3.2) in the case of analytic $V$.    In general it can be realized as a boundary value of the Stieltjes transform of the equilbrium measure and its derivatives.  In the case of sufficiently regular equilibrium measures one can show that it is a sum of delta functions at the spectral edges and a continuous density in the interior.

The proof of Theorem \ref{thm: beta} is now nearly identical to the proof of Theorem \ref{thm1}, replacing Propostion \ref{prop:stein} with Theorem \ref{thm: beta-clt}.   We apply the homogenization result Theorem \ref{thm:betahomog} to the re-scaled $\beta$-ensembles,
\beq
\frac{ \rhoV ( \gV_i ) }{\rhosc (0)} ( x_j - \gV_i ), \qquad \frac{ \rho^{(U)} ( \gU_i ) }{ \rhosc (0) } ( y_j - \gU_i)
\eeq
where $x_i$ is the $\beta$-ensemble under consideration and $y_j$ is a Gaussian $\beta$-ensemble.

 The only thing that needs to be checked is that the variance of $V ( \hat{\G})$ does not depend on the parameters $A, B$ (up to $o ( \log(N) )$ terms).  Using the fact that $\hat{\G'} \neq 0$ only near $\gV_i$, a short calculation using the bounds \eqref{eqn:Gbds} (which also hold for $\hat{\G}$), one sees that for any $\kappa >0$,
\beq
V ( \hat{\G} ) = \frac{1}{2 \beta \pi^2} \int_{|y-\gV_i | \le \kappa } \int_{ |x- \gV_i | \le \kappa } \left( \frac{\hat{\G} (x) - \hat{\G} (y) }{ x-y } \right)^2 \d x \d y + \O (1).
\eeq
The double integral does not depend on $A$ and $B$, and so the proof proceeds as before.  


\subsection{Proof of Theorems \ref{thm:evalues} and Theorem \ref{thm:betaexpect1}}
We note the following elementary definite integrals, which correspond to the terms on the right side of \eqref{eqn:Etr}.
\begin{prop}\label{eqn: definite}
 Let $\gamma\in (-2,2)$ and $f$ be the indicator function:
  \[f(x) = \1_{ \{ x \le \gamma \}}.\]
Then,
\beq
\frac{1}{2 \pi} \int_{-2}^2 f(x) \frac{x^4-4x^2+2}{ \sqrt{4-x^2}} \d x = \frac{\sqrt{4-\gamma^2}}{8 \pi} (2\gamma-\gamma^3),
\eeq
\beq
  \frac{1}{ 2 \pi } \int_{-2}^2 f (x) \frac{2-x^2}{ \sqrt{4-x^2}} \d x= \frac{\sqrt{4-\gamma^2}}{4 \pi} \gamma,
  \eeq
  \beq
  \frac{1}{ 2 \pi } \int_{-2}^2 f (x) \frac{x}{ \sqrt{4-x^2}} \d x= -\frac{\sqrt{4-\gamma^2}}{2 \pi}.
  \eeq
  \qed
  \end{prop}
  
  In order to prove the expansion \eqref{eqn:expect}, we will use the following estimate, which is a consequence of gap universality \cite{EY}.  This estimate states that there is a constant $\mathfrak{a} >0$ so that for all indices $i, j$ with $\kappa N \le i, j \le (1- \kappa)N$ where $\kappa >0$,
  \beq \label{eqn:gu}
 \left|  \ee_W [ \rhosc ( \gamma_i ) ( \lambda_{i+1} - \lambda_i ) ]  - \ee_{GOE}[ \rhosc ( \gamma_j) ( \lambda_{j+1} - \lambda_j ) \right] \le \frac{C}{N^{1+ \mathfrak{a}}}
  \eeq
  where the first expectation is with respect to a Wigner matrix and the second is with respect to a GOE matrix.  Taking now a fixed Wigner matrix with eigenvalues $\lambda_i$, we fix a small $\omega >0$ and write
  \beq
  \lambda_{i_0} = \frac{1}{2 N^\omega+1} \sum_{ |i - i_0 | \le N^{\omega}} \lambda_i + \frac{1}{ 2 N^\omega+1}\sum_{ |i - i_0 | \le N^{\omega}}  \lambda_{i_0} - \lambda_i.
  \eeq
We compare the expectation of $\lambda_{i_0} - \lambda_{i_0+k}$ to $\lambda_{i_0} - \lambda_{i_0-k}$.  We can rewrite each of these quantities as a telescoping sum of gaps of consecutive eigenvalues and apply \eqref{eqn:gu}.  As a consequence, using the smoothness of the semicircle distribution we see that by taking $\omega$ small enough, depending on $\mathfrak{a}$, that
\beq
\left| \ee \frac{1}{ 2 N^\omega+1}\sum_{ |i - i_0 | \le N^{\omega}}  \lambda_{i_0} - \lambda_i \right| \le \frac{C}{N^{1 + \mathfrak{a}/2}}.
\eeq
Hence, it suffices to take the expectation of a local average of eigenvalues.  Fix a small $\eps_1 < \omega /10$ and consider the following smooth function $\varphi$.  We let 
\begin{align}
\varphi (x) = x-\gamma_{i_0}, \qquad & \gamma_{i_0-(N^{\omega} +N^{\eps_1})} \le x \le \gamma_{i_0+N^{\omega}+N^{\eps_1}} \notag\\
\varphi(x) = \gamma_{i_0+N^{\omega}+2N^{\eps_1}} -\gamma_{i_0}, \qquad & x \ge \gamma_{i_0+N^{\omega}+3 N^{\eps_1}} \notag\\
\varphi (x) = \gamma_{i_0-(N^{\omega}-2N^{\eps_1})} - \gamma_{i_0}, \qquad & x \le \gamma_{i_0-(N^{\omega}+3N^{\eps_1})}.
\end{align}
We let $\varphi$ smoothly interpolate between these values so that $|\varphi'(x) | \le C$, and $|\varphi^{(k)} (x) | \le C N^{(k-1)(1-\eps_1)}$ for $k=2, 3$.  Note that if 
\beq
f(x) = \frac{N}{2 N^{\omega}+1} \varphi (x)
\eeq
then $|| f(x) ||_1 \le C$, $||f'(x) ||_1 \le C$ and $||f''(x) ||_1 \le N^{\omega}/N$ and so we can apply Lemma \ref{lem:var} to $f$. 

Let $I$ denote the interval
\beq
I = [ \gamma_{i_0-N^{\omega}-1}, \gamma_{i_0+N^{\omega}}].
\eeq
By rigidity we have with overwhelming probability that
\begin{align}
\frac{1}{ 2 N^{\omega}+1} \left|  \sum_{ |i-i_0| > N^{\omega} } \varphi ( \lambda_i ) - N \int_{ x \notin I } \varphi (x) \rhosc (x) \d x \right| \le \frac{C N^{2 \eps_1}}{N^{1+\omega}} \le \frac{C}{N^{1+\omega/2}}
\end{align}
by our choice of $\eps_1$.   Using the fact that
\beq
\gamma_i - \gamma_j = \frac{i-j}{N \rhosc ( \gamma_i ) } + \O \left( \frac{(i-j)^2}{N^2} \right)
\eeq
one sees that
\beq
\frac{N}{2 N^{\omega}+1} \int_I \varphi (x) \rhosc (x) \d x= \frac{N}{2 N^{\omega}+1} \int_I (x - \gamma_{i_0} )\rhosc (x) \d x= \frac{-1}{2 \rho ( \gamma_{i_0} ) N } + \O ( N^{2 \omega-2} ) 
\eeq
Hence, we see that for some $c>0$ we have
\beq
N \ee[ \lambda_{i_0}  - \gamma_{i_0}] =  \sum_{i} \ee[ f ( \lambda_i ) ] - N \int f(x) \rhosc (x) \d x -  \frac{1}{2 \rhosc ( \gamma_{i_0} ) } + \O (N^{-c})
\eeq
Note that,
\beq
f(x) \to \frac{1}{2 \rhosc ( \gamma_{i_0} ) } \left( \1_{ \{x \ge \gamma_{i_0} \} } - \1_{ \{ x \le \gamma_{i_0} \}}\right).
\eeq
and so we see from Lemma \ref{lem:var} that
\beq
\sum_i \ee [ f ( \lambda_i ) ] - N \int f(x) \rhosc (x) \d x = \frac{1}{2 \pi \rhosc ( \gamma_i ) } \arcsin \left( \frac{ \gamma_i}{2} \right) - \frac{1-a_2}{2 } \gamma_i + \frac{s_4}{4} (  \gamma^3_i -2  \gamma_i ) + o(1)
\eeq
This yields the claim.  

%

To obtain the result for the variance, we rearrange \eqref{eqn:hzeta} into:
\begin{align*}
x_i(t)-\mathbb{E}[x_i(t)]+\frac{1}{N}\zeta_{y}&=y_i(t)-\mathbb{E}[y_i(t)]\\
&+\frac{1}{N}(\zeta_{x}-\mathbb{E}[\zeta_{x}])\\
&+\frac{1}{N}\mathbb{E}[\zeta_{x}]+\mathbb{E}[y_i(t)]-\mathbb{E}[x_i(t)]+\O(N^{-1-c}).
\end{align*}
Squaring and using the independent between $x_i(t)$ and $\zeta_y$, and $y_i(t)$ and $\zeta_{x}$, we have
\[\mathrm{Var}(x_i(t))+\frac{1}{N^2}\mathrm{Var}(\zeta_{y})+\frac{1}{N^2}(\mathbb{E}[\zeta_{y}])^2=\mathrm{Var}(y_i(t))+\frac{1}{N^2}\mathrm{Var}(\zeta_{x})+(\mathbb{E}[\frac{1}{N}\zeta_{x}-x_i(t)+y_i(t)])^2+\O(N^{-2-2c}).\]

By \eqref{eqn: rep}, we have
\[\frac{1}{N^2}(\mathbb{E}[\zeta_{y}])^2=\big(\mathbb{E}[\frac{1}{N}\zeta_{x}-x(t)+y(t)]\big)^2+\O(N^{-2-2c}),\]
so
\begin{align*}
&\mathrm{Var}(x_i(t))=\mathrm{Var}(y_i(t))+\frac{1}{N^2}(\mathrm{Var}(\zeta_{x})-\mathrm{Var}(\zeta_{\mathbf{y}}))+\O(N^{-2-2c})\\
=&\mathrm{Var}(y_i(t))+\frac{a_2-1}{4 \pi^2N^2} \left( \int_{-2}^2 \G(x) \frac{x}{\sqrt{4-x^2} }\,\mathrm{d}x \right)^2 + \frac{s_4}{2 \pi^2N^2} \left( \int_{-2}^2 \G(x) \frac{ 2-x^2}{\sqrt{4-x^2} } \,\mathrm{d}x \right)^2+\O(N^{-2-2c})
\end{align*}
Using \eqref{eqn: approximate-indicator} and Proposition \ref{eqn: definite}, we have
\begin{align*}
&\frac{a_2-1}{4 \pi^2} \left( \int_{-2}^2 \G(x) \frac{x}{\sqrt{4-x^2} }\,\mathrm{d}x \right)^2 + \frac{s_4}{2 \pi^2} \left( \int_{-2}^2 \G(x) \frac{ 2-x^2}{\sqrt{4-x^2} } \,\mathrm{d}x \right)^2\\
=~&(a_2-1)\gamma_i^2 + \frac{s_4}{8}+o(1).
\end{align*}
The gap universality result of \cite{EY} also holds for $\beta$-ensembles, and so the proof of Theorem \ref{thm:betaexpect1} is similar to the Wigner matrix case. \qed

\subsection{Proof of Theorem \ref{thm: BPZ}}

First we consider the case \eqref{eqn: BPZ-1}.  We will smooth out the indicator function $\1_{ \{ x \le u \}}$.  Let $ \chi (x)$ be a smooth function so that $\chi (x) = 1$ for $x \le u$ and $\chi (x) = 0$ for $ x \ge u + N^{\omega-1}$ where $0 < \omega < 1/10$, and $| \chi^{(k)} (x) | \le C N^{k( \omega-1)}$, $k=1, 2, 3$.  Let $i_0$ be the index of the classical eigenvalue closest to $u$.  By rigidity for any $\eps >0$, we have with overwhelming probability
\beq
\left| \sum_i f ( \lambda_i ) \1_{ \{ \lambda_i \le u \} } - \sum_i f ( \lambda_i ) \chi ( \lambda_i )  \right| \le \sum_{  |i-i_0| \le N^{\omega+\eps}} |f ( \lambda_i ) | \le  C \frac{N^{2 \omega+2 \eps}}{N},
\eeq
where we used $f(u) = 0$, $|f'| \le C$ and $|\lambda_i - u | \le C N^{\omega+\eps}/N$ for $i$ appearing in the rightmost sum.  It then suffices to apply Proposition \ref{prop:stein} to the function $f(x) \chi (x)$.


The proof of \eqref{eqn: BPZ-2} is similar.  First, by subtracting a constant we may assume that $f ( \gamma_k) = 0$.  Let $\chi(x)$ be the indicator function as above, but for $u = \gamma_k$.  Then, for any $\eps >0$, we have with overwhelming probability,
\beq
\left| \sum_{i=1}^k f ( \lambda_i) - \sum_{i=1}^N f ( \lambda_i ) \chi ( \lambda_i ) \right| \le \sum_{ |i - k| \le N^{ \omega+\eps}} | f ( \lambda_i ) | \le C \frac{N^{2 \omega+2 \eps}}{N}.
\eeq
The result follows from applying Proposition \ref{prop:stein} to $f(x) \chi (x)$.


\bibliography{meso_bib}{}
\bibliographystyle{abbrv}

\end{document}